\documentclass[11pt]{article}
\usepackage{amssymb}
\usepackage{bbm}
\newcommand{\Boxneu}{\square} 
\newcommand{\nc}[2]{\newcommand{#1}{#2}}

\nc{\bsa}{\begin{satz}}
\nc{\bpr}{\begin{pr}}
\nc{\bth}{\begin{them}}
\nc{\ble}{\begin{lem}}
\nc{\bco}{\begin{corollary}}
\nc{\bre}{\begin{remark}}
\nc{\bex}{\begin{example}}
\nc{\bde}{\begin{de}}
\nc{\ede}{\end{de}}
\nc{\esa}{\end{satz}}
\nc{\epr}{\end{pr}}
\nc{\ethe}{\end{them}}
\nc{\ele}{\end{lem}}
\nc{\eco}{\end{corollary}}
\nc{\ere}{\hfill\mbox{$\Diamond$}\end{remark}}
\nc{\eex}{\end{example}}
\nc{\epf}{\hfill\mbox{$\Boxneu$}}
\nc{\beq}{\begin{equation}}
\nc{\eeq}{\end{equation}}
\nc{\ot}{\otimes}
\nc{\lra}{\longrightarrow}
\nc{\ci}{\circ}

\begin{document}
\title{Gauge transformations on locally trivial quantum principal fibre
bundles}
\author{Dirk Calow$^1$\thanks{supported by Deutsche
Forschungsgemeinschaft, e-mail Dirk.Calow@itp.uni-leipzig.de}
\ and Rainer Matthes$^{1,2}$ \thanks{supported by S\"achsisches
Staatsministerium f\"ur Wissenschaft und Kunst,\newline\hspace*{.5cm} e-mail
Rainer.Matthes@itp.uni-leipzig.de or rmatthes@mis.mpg.de}\\[.5cm]
\normalsize $^1$Institut f\"ur Theoretische Physik
\normalsize der Universit\"at Leipzig\\
\normalsize Augustusplatz 10/11,
\normalsize D-04109 Leipzig,
\normalsize Germany
\\[.5cm]
\normalsize $^2$Max-Planck-Institut f\"ur Mathematik
\normalsize in den Naturwissenschaften\\
\normalsize Inselstra{\ss}e 22-26,
\normalsize D-04103 Leipzig,
\normalsize Germany}
\date{}
\maketitle
\begin{abstract}
\newtheorem{pr}{Proposition}
\newtheorem{de}{Definition}
\newtheorem{lem}{Lemma}
\newtheorem{them}{Theorem}
We consider in this paper gauge 
transformations on locally trivial quantum principal fibre bundles (QPFB). 
If ${\cal P}$, $B$, $H$, are the algebras of the total space, the base space,
the structure group of the bundle, left (right) gauge transformations are
defined as
 automorphisms of the left (right) $B$-module ${\cal P}$ which are
adapted to the coaction of the
 Hopf algebra $H$ and to the covering related
to the local trivializations. 
 Connections on the QPFB are in general not
transformed into connections.
 For covariant connections, there are analogues
of the classical formulas
 relating the connection and its gauge transform. 
\end{abstract}
This paper is a follow-up of 
\cite{cama}--\cite{cama2}. We freely use the results of these
papers. In \cite{Kon}, gauge transformations are defined as ``vertical
automorphisms'' of the bundle. We start here with the more general definition
of \cite{brma} (given there only for trivial bundles) and \cite{haj} of gauge
transformations as convolution invertible linear maps from the Hopf algebra
$H$ to the  total space algebra ${\cal P}$ (taken now in a different
context) and consider the action of such  transformations on differential
geometric objects (connections and covariant derivatives).

We define left (right) gauge transformations on a locally trivial QPFB in the
sense of  \cite{Kon} and \cite{cama1} as automorphisms of the left (right)
$B$-module ${\cal P}$ which are compatible with the right coaction on the
total space and with the local trivializations. On a trivial bundle,
such gauge transformations correspond to convolution invertible (left case)
and twisted convolution invertible (right case) linear maps from
$H$ to $B$. On a locally trivial QPFB, the corresponding objects
are linear maps from $H$ to ${\cal P}$ with these properties. 
If the structure group is a compact quantum group, gauge transformations
can be characterized locally by elements of the algebras defined locally
by the covering of the basis. Gauge transformations can also be defined
on locally trivial quantum vector bundles (QVB) as module automorphism
respecting the local trivializations.  If a QVB is associated to some
QPFB, every gauge transformation of the QPFB induces a gauge transformation
on the QVB.

Every gauge transformation on a QPFB induces a module isomorphism of the
module of horizontal forms belonging to a differential structure of ${\cal P}$
(the module structure being with respect to the maximal embeddable LC
differential algebra related to the differential structure). 
It follows that the set of covariant derivatives is invariant under
gauge transformations. For the set of connections, this is true at least
in the following two cases: If the differential structure on ${\cal P}$
is defined using the universal calculus on $H$, and if the differential
structure on ${\cal P}$ uses a bicovariant calculus on $H$, restricting
at the same time to  ``vertical automorphisms''. In any case, the local
connection forms and the curvature forms are transformed in the typical
(inhomogeneous respectively homogeneous) way.
We also give an example of a gauge transformation on a $U(1)$ bundle over a
gluing of $SU_\nu(2)$ with some ``quantum cylinder'', which is not an algebra
automorphism.

\section{\label{gauge} Gauge transformations}
In \cite{Kon} gauge
transformation on a QPFB $\cal P$ are defined as automorphisms $\alpha :
{\cal P} \longrightarrow {\cal P}$ fulfilling the conditions
\begin{eqnarray} (\alpha \otimes id) \circ \Delta_{\cal P} &=&\Delta_{\cal
P} \circ \alpha \\ \alpha \circ \iota &=& \iota.
\end{eqnarray}
On trivial QPFB gauge transformations are in one to one
correspondence with homomorphisms $\tau_{\alpha} : H \longrightarrow B$
fulfilling
\begin{eqnarray*}
\tau_{\alpha}(1)&=& 1 \\ \tau_{\alpha}
(h)a &=& a\tau_{\alpha}(h), \;\forall h
\in H\; \forall a \in B 
\end{eqnarray*}
such that 
\[
\alpha(a \otimes h)=\sum 
a\tau_{\alpha}(h_1) \otimes h_2. \]
To get nonclassical gauge transformations one needs a
more general definition.
First, let us define gauge transformation on trivial QPFB.
\begin{de}  Let $B \otimes H$ be a trivial QPFB. A left (right) gauge 
transformation on $B \otimes H$
is a left (right) $(B \otimes 1)$-module isomorphism 
$\alpha : B \otimes H \longrightarrow B \otimes 
H$ satisfying \begin{eqnarray} \label{alpha}(\alpha\otimes id) \circ (id
\otimes 
\Delta)&=&(id \otimes \Delta) \circ \alpha \\ \alpha \circ (id \otimes 
1)&=&(id \otimes 1). \end{eqnarray}
\end{de}
Remark: Obviously, for every gauge transformation $\alpha$, the
inverse 
$\alpha^{-1}$ is also a gauge transformation.

In the sequel, our standard notation will be $\alpha_l$ for left and 
$\alpha_r$ for right gauge transformations. The notation $\alpha_{l,r}$
will be used if something is true for both left and right gauge 
transformations. 
\begin{pr} Left (right) gauge transformations $\alpha_{l,r}$ on a trivial
bundle $B\otimes H$ are in one to
one 
correspondence 
to linear maps $\tau_{\alpha_{l,r}} : H \longrightarrow B$ satisfying 
\begin{eqnarray*} \tau_{\alpha_{l,r}}(1)&=& 1 \\
\sum \tau_{\alpha^{-1}_l}(h_1) \tau_{\alpha_l}(h_2) &=&\sum
\tau_{\alpha_l}(h_1) 
\tau_{\alpha^{-1}_l}(h_2)= \varepsilon(h)1 \\
\sum \tau_{\alpha^{-1}_r}(h_2) 
\tau_{\alpha_r}(h_1)&=&\sum \tau_{\alpha_r}(h_2)\tau_{\alpha^{-1}_r}(h_1)=
\varepsilon(h)1 \end{eqnarray*}
such that \begin{eqnarray} \alpha_l(a \otimes h)&=& \sum a
\tau_{\alpha_l}(h_1) 
\otimes h_2 \\ \alpha_r(a \otimes h)&=& \sum \tau_{\alpha_r}(h_1)a \otimes
h_2 .
\end{eqnarray} \end{pr}
Proof: 
 We define 
\begin{equation} \tau_{\alpha_{l,r}}(h) =(id \otimes 
 \varepsilon) \circ \alpha_{l,r}(1 \otimes h).
 \end{equation} 
 Because of
formula (\ref{alpha}) and  
 $(id \otimes \varepsilon \otimes id) \circ (id \otimes \Delta)=id$ 
 we obtain for left gauge transformations \begin{eqnarray*}
 \alpha_l(a \otimes h) &=&(id \otimes \varepsilon \otimes id) \circ (id
\otimes 
 \Delta)    \circ \alpha_l(a 
 \otimes h) \\
 &=& (a \otimes 1)(id \otimes \varepsilon \otimes id) \circ 
 (\alpha_l \otimes id) \circ (id \otimes \Delta)(1 \otimes h) \\ &=& (a 
 \otimes 
 1)\sum (\tau_{\alpha_l}(h_1) \otimes h_2) \\
 &=& \sum (a \tau_{\alpha_l}(h_1) \otimes h_2) \end{eqnarray*}
 and for right gauge transformations 
 \begin{eqnarray*}
 \alpha_r(a \otimes h) &=&(id \otimes \varepsilon \otimes id) \circ (id
\otimes 
 \Delta)    
 \circ \alpha_r(a 
 \otimes h) \\
 &=&(id \otimes \varepsilon \otimes id) \circ 
 (\alpha_r \otimes id) \circ (id \otimes \Delta)(1 \otimes h)(a \otimes
1)
  \\ &=& \sum (\tau_{\alpha_r}(h_1) \otimes h_2)(a \otimes 1) \\
  &=& \sum (\tau_{\alpha_l}(h_1)a \otimes h_2). \end{eqnarray*}
  Since the gauge transformations are 
  left (right) $(B \otimes 1)$-module isomorphisms, the 
properties claimed for 
  $\tau_{\alpha_{l,r}}$ are easily verified.
  We leave the other direction of the proof to the reader. \hfill$\Boxneu$ 
  \begin{de} Let $\cal P$ be a locally trivial QPFB with local
trivializations
$\chi_i:{\cal P}\lra B_i\ot H$.\\
 A left (right) gauge transformation on $\cal P$ is a left
$\iota(B)$-module 
isomorphism 
$\alpha_l,r: {\cal P} \longrightarrow {\cal P}$ 
such that there exists a family $(\alpha_{{l,r}_i})_{i \in I}$ 
  of left (right) gauge transformation on the trivializations $B_i \otimes
H$
satisfying 
  \begin{equation} \label{alphi} \chi_i \circ \alpha_{l,r} =
\alpha_{{l,r}_i} 
\circ\chi_i. \end{equation}
\end{de}
  \begin{pr} A left (right) gauge transformation fulfills \begin{equation}
  \Delta_{\cal P} \circ \alpha_{l,r} = (\alpha_{l,r} \otimes id) \circ 
  \Delta_{\cal P}. \label{gtr}\end{equation} \end{pr}
  Proof: Using the definitions one calculates
  \begin{eqnarray*}  (\chi_i \circ \alpha_{l,r} \otimes id)\circ
\Delta_{\cal P} 
  &=& 
  (\alpha_{{l,r}_i} \circ \chi_i \otimes id) \circ \Delta_{\cal P} \\
  &=& ((\alpha_{{l,r}_i} \otimes id) \circ (id \otimes \Delta) \circ \chi_i
\\
  &=& (id \otimes \Delta) \circ \alpha_{{l,r}_i} \circ \chi_i \\
  &=& (id \otimes \Delta) \circ \chi_i \circ \alpha_{l,r} \\
  &=& (\chi_i \otimes id) \circ \Delta_{\cal P} \circ \alpha_{l,r}. 
  \end{eqnarray*}
 (\ref{gtr}) follows from $\bigcap_i ker \chi_i=0$. \hfill$\Boxneu$
  \\\\Remark: Left (right) gauge transformations $\alpha_{l,r}$ can be 
  equivalently
  defined as left (right) $\iota(B)$ module isomorphisms $\alpha_{l,r}:
{\cal P} 
  \longrightarrow {\cal P}$ fulfilling \begin{eqnarray*}			   
  \Delta_{\cal P} \circ \alpha_{l,r} &=& (\alpha_{l,r} \otimes id) \circ 
  \Delta_{\cal P} \\ \alpha_{l,r} \circ \iota &=& \iota \\
  \alpha(ker \chi_i)&=& ker \chi_i;~~\forall i \in I. \end{eqnarray*}
  \begin{pr} The set $G_{l,r}$ of all left (right) transformations is a
group 
  with 
  the composition of maps as group multiplication. \hfill$\Boxneu$ \end{pr}
  \begin{pr} \label{gauone} Left (right) gauge transformations
$\alpha_{l,r}$ on 
  a locally 
  trivial QPFB $\cal P$
  are in one to one correspondence to linear maps $g_{\alpha_{l,r}} : H 
  \longrightarrow \cal P$ satisfying \begin{eqnarray}
  \label{gauone1} g_{\alpha_{l,r}}(1) &=& 1\\
  \label{gauone2} \Delta_{\cal P}(g_{\alpha_l}(h))&=& \sum
g_{\alpha_l}(h_2) 
  \otimes S(h_1)h_3 \\
  \label{gauone3} \Delta_{\cal P}(g_{\alpha_r}(h) &=& \sum
g_{\alpha_r}(h_2) 
  \otimes h_3 
  S^{-1}(h_1)\\ 
  \label{gauone4} \sum g_{\alpha_l}(h_1) 
  g_{\alpha^{-1}_l}(h_2)&=&\sum g_{\alpha^{-1}_l}(h_1)g_{\alpha_l}(h_2) 
  =\varepsilon(h)1 \\
  \label{gauone5} \sum g_{\alpha_r}(h_2) g_{\alpha^{-1}_r}(h_1) &=& \sum 
  g_{\alpha^{-1}_r}(h_2) 
  g_{\alpha_r}(h_1) = \varepsilon(h) 1. \end{eqnarray}
The correspondence is given by
 \begin{eqnarray} \alpha_l(f) &=& \sum f_0g_{\alpha_l}(f_1) \\
  \alpha_r(f) &=& \sum g_{\alpha_r}(f_1)f_0.  \end{eqnarray} \end{pr}
  Proof: We will give the proof only for left gauge transformations 
  because it works for right gauge transformations with the same arguments.
\\Assume that there is given a linear map $g : H \longrightarrow \cal P$
with 
  the properties $\Delta_{\cal P}(g(h))=\sum g(h_2) \otimes S(h_1)h_2$ and 
  $g(1)=1$ such 
  that there exists a linear map $g^{-1}: H \longrightarrow \cal P$
fullfilling
  $\sum g^{-1}(h_1)g(h_2)=\sum g(h_1)g^{-1}(h_2)=\varepsilon(h)1$.
 It is easy to verify that the linear map $\alpha : {\cal P}
\longrightarrow 
  {\cal P}$ defined by \[ \alpha(f):=\sum f_0g(f_1) \] is a left gauge 
  transformation.
  \\The proof of the other direction is more complicated.
  By definition for a given left gauge transformation 
  $\alpha_l$ there exists a 
  family $(\alpha_{l_i})_{i \in I}$ of left gauge transformations on the 
  trivializations 
  $B_i \otimes H$. Since the linear maps $\alpha_{l_i}$ are left $(B_i
\otimes 
  1)$-module isomorphisms, there exist left $(B_{ij} \otimes 1)$ module
isomorphisms $\alpha^j_{l_i} : B_{ij} \otimes H \longrightarrow B_{ij} 
  \otimes H$ satisfying \begin{equation} \label{alphaij} \alpha^j_{l_i} 
  \circ (\pi^i_j \otimes id) =(\pi^i_j 
  \otimes id) \circ \alpha_{l_i}. \end{equation} 
These $\alpha^j_{l_i}$ satisfy the identity
  \begin{equation} \label{gi} \alpha^j_{l_i} = \phi_{ij} \circ
\alpha^i_{j_l} 
  \circ \phi_{ji}, \label{aijl}\end{equation}
where $\phi_{ij}$ are the isomorphisms $\phi_{ij}: B_{ij} \otimes H \longrightarrow B_{ij} \otimes 
H $ induced from the transition functions $\tau_{ij}$ of the bundle ${\cal P}$ 
(see \cite{cama2}). (\ref{aijl}) is proved as follows.
 Let $f \in {\cal P}$. We know, \begin{equation} \label{shit} 
  (\pi^i_j \otimes id) \circ \chi_i(f) = \phi_{ij} \circ (\pi^j_i \otimes
id) \circ \chi_j(f), \end{equation} therefore \[
  (\pi^i_j \otimes id) \circ \chi_i(\alpha_l(f)) = \phi_{ij} \circ 
  (\pi^j_i \otimes id) \circ \chi_j(\alpha_l(f)). \] 
With (\ref{alphi}) and (\ref{alphaij}) follow the equations
 \begin{eqnarray}  (\pi^i_j \otimes id) \circ \alpha_{l_i} \circ \chi_i(f) 
 &=& \phi_{ij} \circ 
  (\pi^j_i \otimes id) \circ \alpha_{j_l} \circ \chi_j(f) \\ \label{shit1}
  \alpha^j_{l_i} \circ (\pi^i_j \otimes id) \circ \chi_i(f) &=& \phi_{ij}
\circ \alpha^i_{j_l} \circ (\pi^j_i \otimes id) \circ \chi_j(f). \end{eqnarray}
Inserting (\ref{shit}) in (\ref{shit1}) one obtains 
  \[ \alpha^j_{l_i} \circ \phi_{ij} =\phi_{ij} \circ \alpha^i_{j_l}, \]
which proves (\ref{aijl}). Because of this formula, the linear maps 
  $\tau_{\alpha_{l_i}}: H \longrightarrow B_i$ corresponding to the 
  $\alpha_{l_i}$ 
   satisfy \begin{equation} \label{gei} \pi^i_j(\tau_{\alpha_{l_i}}(h)) = 
   \sum \tau_{ij}(h_1) 
   \pi^j_i(\tau_{\alpha_{j_l}}(h_2)) \tau_{ji}(h_3). \end{equation}
   Now we define a family of linear maps $g_{\alpha_{l_i}} : H
\longrightarrow B_i 
   \otimes H$ by \[ g_{\alpha_{l_i}}(h) := \sum \tau_{\alpha_{l_i}}(h_2)
\otimes S(h_1)h_3 .\] 
It is easy to see that \[ \alpha_{l_i}(a \otimes
h)=\sum (a \otimes h_1) g_{\alpha_{l_i}}(h_2). \]
   Because of formula (\ref{gei}) the family of linear maps 
   $(g_{\alpha_{l_i}})_{i \in 
   I}$ fulfills \begin{equation} (\pi^i_j \otimes id) (g_{\alpha_{l_i}}(h))
= \phi_{ij} \circ (\pi^j_i \otimes id)(g_{\alpha_{j_l}}(h)),
\end{equation}
   i.e. there exists a unique linear map $g_{\alpha_l} : H \longrightarrow
{\cal P}$  satisfying  $\chi_i(g_{\alpha_l}(h)) = g_{\alpha_{l_i}}(h)$ and 
   $\alpha_l(f) =\sum f_0g_{\alpha_l}(f_1)$. The properties of
$g_{\alpha_l}$ are now easily verified by using
   the properties of the family $(g_{\alpha_{l_i}})_{i \in I}$.
\hfill$\Boxneu$
   \begin{pr} There is a bijection between left and right gauge
transformations. 
   \end{pr}
   Proof: Left and right gauge transformations are related by 
$g_{\alpha_r}:= g_{\alpha_r} \circ S^{-1}$.\hfill$\Boxneu$
\begin{pr} \label{etae} Let $\Gamma({\cal P})$ be a differential
structure on 
   $\cal P$ and let $\alpha_{l,r}$ be a left (right) gauge transformation
on $\cal P$.\\
Then the formulas \begin{eqnarray} \label{curv3} \alpha_l(\gamma)&:=& 
\sum \gamma_0 g_{\alpha_l} (\gamma_1)\\ \label{curv4}
   \alpha_r(\gamma)&:=& \sum g_{\alpha_r}(\gamma_1) \gamma_0 \end{eqnarray}
define left (right) $\Gamma_m(B)$-module isomorphisms
   $\alpha_{l,r} : hor \Gamma_c({\cal P}) \longrightarrow hor
\Gamma_c({\cal P})$ . \end{pr}
($\Gamma_m(B)$ is the maximal embeddable LC-differential algebra induced from 
$\Gamma_c({\cal P})$, see \cite{cama2}.) We leave the proof to the reader.
\\The concept of gauge transformations can be carried over to associated
vector bundles. 
In general one can define gauge transformations on a locally trivial QVB
$E$ as follows: 
\begin{de} Let $((E,B,\kappa),V,(\zeta_i,J_i)_{i \in I})$ be a
locally trivial QVB. A gauge transformation on $E$ is a automorphism $\eta:
E\longrightarrow E$ with the properties \begin{eqnarray} \kappa(a) \circ
\eta &=& \eta \circ \kappa(a),\; \forall a \in B \\
\eta(ker \zeta_i) &=& ker \zeta_i. \end{eqnarray} \end{de}
The last condition of this definition has the consequence that there are
gauge transformations $\eta_i : B_i \otimes V \longrightarrow B_i
\otimes V$ and 
   $\eta^i_{ij} : B_{ij} \otimes V \longrightarrow B_{ij} \otimes V$ such
that \begin{eqnarray*} \eta_i \circ \zeta_i &=& \zeta_i \circ \eta \\
   \eta^i_{ij} \circ (\pi^i_j \otimes id) &=& (\pi^i_j \otimes id) \circ
\eta_i. \end{eqnarray*}

\begin{pr}
Gauge transformations on a locally trivial QVB $E$ are in
one-to-one correspondence with families of gauge transformations 
$\eta_i : B_i \otimes V \longrightarrow B_i \otimes V$ 
satisfying \begin{equation} \eta^i_{ij} = \phi_{{ij}_E} \circ 
\eta^j_{ij} \circ \phi_{{ji}_E}. \end{equation} 
\end{pr}
   We omit the proof because it is quite analogous to the proof 
 of Proposition 4 of \cite{cama2}.

\begin{pr} 
Let $E({\cal P},F)$ be an associated vector bundle. Every
left gauge transformation on $\cal P$ determines a gauge transformation 
on $E({\cal P},F)$. 
\end{pr}
Proof: Let $\alpha_l$ be a left gauge transformation on $\cal P$. 
The linear map $\eta_{\alpha_l} : E({\cal P},F)
\longrightarrow 
   E({\cal P},F)$ 
  defined by 
\[ 
\eta_{\alpha_l}:= \alpha_l \otimes id 
\]
 is seen to be a
   gauge transformation on $E({\cal P},F)$. \hfill$\Boxneu$
   \\\\ Let $\Gamma_m(B)$ be the maximal embeddable LC-differential algebra over $B$
induced from the differential structure on $\cal P$, and let 
$E_{\Gamma}({\cal P},F)$ be
the locally trivial QVB constructed in terms of $E({\cal P},F)$ and 
   $\Gamma_m(B)$. Because of Proposition \ref{etae}, one can extend e
 very gauge transformation 
on $E({\cal P},F)$ determined by a gauge transformation on $\cal P$ to a
module automorphism $\eta_{\alpha}$ of $E_{\Gamma}({\cal P},F)$ by 
\[ \eta_{\alpha_l}:= \epsilon_{\Gamma} \circ (\alpha_l \otimes id) \circ 
\epsilon^{-1}_{\Gamma}.\]
We end up this section with some remarks about the general structure
of gauge transformations in the case when the structure group of the locally
trivial QPFB is a compact quantum group. 
   \\For the algebra $P(G)$ of polynomial functions over such a
compact quantum group $G$ one can 
   construct the following linear basis ( see \cite{wor3}, \cite{korn} 
   and \cite{dipl}). 
   Let $M$ be the set of all irreducible unitary matrix co-representations
of $G$. \\
(A unitary matrix co-representation is defined by an $P(G)$-valued
$N \times N$-matrix $(u_{ij})_{i,j=1,2,...,N}$ with
$\Delta(u_{ij})=\sum_k 
   u_{ik} \otimes u_{kl}$ and $u^*_{ij}=S(u_{ji})$.)
   \\ Two co-representations $\rho$ and $\sigma$ are equivalent if there
exists an intertwining operator. This defines an 
   equivalence relation $\sim$ in $M$. Now one can select in every class 
   $\alpha=M/\sim$ a matrix co-representation 
   $(u^{\alpha}_{ij})_{i,j=1,2,...N^{\alpha}}$. It is proved in
\cite{wor3}, 
   \cite{korn} and {\cite{dipl} that the set of elements
   $(u^{\alpha}_{ij})_{\alpha \in M/\sim;\,i,j=1,2,...,N^{\alpha}}$ is a 
   linear 
   basis of $P(G)$.
   This leads us to the following conlusion. 
\bpr
Let $\cal P$ be a locally
trivial 
   QPFB 
   where the structure group is a compact quantum group. The set of
left (right) gauge 
transformations is in one to one correspondence with sets of invertible 
$B_i$-valued matrices $(b^{\alpha}_{i_{kl}})_{k,l=1,2,...N^{\alpha} 
\;\alpha 
\in M/\sim}$ satisfying 
\begin{equation} 
\label{invma} 
   \pi^i_j(b^{\alpha}_{i_{kl}})=\sum^{N^{\alpha}}_{m,n}   
\tau_{ij}(u^{\alpha}_{km})\pi^j_i(b^{\alpha}_{j_{mn}})\tau_{ji}(u^{\alpha}_{ml}). 
\end{equation}
\epr
 More precisely (see also formula (\ref{gei})), one can construct a gauge
transformation by mapping 
 every matrix $u^{\alpha}_{ij}$ to a set of invertible $B_i$-valued 
 $N^{\alpha} \times N^{\alpha}$-matrices satisfying (\ref{invma}). Doing
this for each $\alpha$ one obtains linear maps $\tau_i : P(G)
\longrightarrow  B_i$ 
   defined by \[ \tau_i(u^{\alpha}_{kl})=b^{\alpha}_{i_{kl}}. \] 
Since the matrices $(b^{\alpha}_{i_{kl}})$ are invertible the $\tau_i$ are 
   convolution invertible. Since formula (\ref{invma}) is fulfilled for
each $\alpha$ the $\tau_i$ satisfy also (\ref{gei}) and it follows that
they determine a gauge transformation.	
\section{Gauge transformations and connections}
Since we have extended gauge transformations only to the subalgebra of 
horizontal forms $hor \Gamma_c({\cal P})$ it is not possible to transform a
connection by a gauge transformation analogous to the classical case by 
transforming the connection form $\omega_{l,r}$ . 
What we can transform is the covariant derivation $D_{l,r}:=hor \circ d$ 
corresponding to the connection.
\\
Let $D_{l,r}$ be the left (right) covariant derivative corresponding to a
left 
(right) connection and let $\alpha_{l,r}$ be a left (right) gauge 
transformation. One defines the linear map $D'_{l,r}:  hor \Gamma_c({\cal P}) 
\longrightarrow hor \Gamma_c({\cal P})$ by 
\[ D'_{l,r} := \alpha_{l,r} \circ D_{l,r} \circ \alpha^{-1}_{l,r}. \]
In the sequel we want to discuss this formula for connections on trivial
QPFB.
We know, that a left (right) connection corresponds to a linear map 
$A_{l,r} : H \longrightarrow \Gamma^1(B)$ satisfying (54), (55),(79)
and (80) of \cite{cama1} respectively. We are interested in the transformed
maps $A'_{l,r}$  belonging to the $D'_{l,r}$.
One calculates \begin{eqnarray*} D'_l(a \otimes h) &=& \alpha_l D_l \circ 
\alpha^{-1}_l(a \otimes h)\\
&=& \alpha_l \circ hor_l \circ d \sum (a \tau_{\alpha^{-1}_l}(h_1) \otimes
h_2) \\
&=& \alpha_l \circ hor_l (\sum ((da) \tau_{\alpha^{-1}_l}(h_1)) \hat{\otimes}
h_2 + \sum a \tau_{\alpha^{-1}_l}(h_1) \hat{\otimes} dh_2) \\
&=& \alpha_l(\sum (d a\tau_{\alpha^{-1}_l}(h_1)) \hat{\otimes} h_2 - 
\sum a \tau_{\alpha^{-1}_l}(h_1)A_l(h_2) \hat{\otimes} h_3) \\
&=& \sum ((da) \tau_{\alpha^{-1}_l}(h_1))\tau_{\alpha_l}(h_2) \hat{\otimes}
h_3 - \sum a \tau_{\alpha^{-1}_l}(h_1)A_l(h_2)\tau_{\alpha_l}(h_3)
\hat{\otimes}
h_4 \\ &=& (da) \hat{\otimes} h -\sum a \tau_{\alpha^{-1}_l}(h_1)d 
\tau_{\alpha_l}(h_2) 
\hat{\otimes} h_3 - \sum a \tau_{\alpha^{-1}_l}(h_1)A_l(h_2) 
\tau_{\alpha_l}(h_3) \hat{\otimes} h_4. \end{eqnarray*}
This shows that the linear map ${A'}_l : H \longrightarrow \Gamma^1(B)$
defined by 
\[ {A'}_l(h):= (id \otimes \varepsilon) \circ D'_l \] has the form
					
\begin{equation} \label{transfa}
A'_l(h)= \sum \tau_{\alpha^{-1}_l}(h_1)A_l(h_2) \tau_{\alpha_l}(h_3) + 
\sum \tau_{\alpha^{-1}_l}(h_1)d \tau_{\alpha_l}(h_2). \end{equation}
The same calculation for right connections leads to \begin{equation}
A'_r(h) = \sum \tau_{\alpha_r}(h_3) A_r(h_2) \tau_{\alpha^{-1}_r}(h_1) - 
\tau_{\alpha_r}(h_2)d \tau_{\alpha^{-1}_r}(h_1). \end{equation}
In general the linear maps $A'_{l,r}$ do not satisfy the conditions 
(79) and (80) of \cite{cama1} respectively, i.e. $D'_{l,r}$ do in
general
not define a connection. Only in the special case when $\Gamma(H)$ is the
universal 
differential algebra, which means every linear map $A :H \longrightarrow 
\Omega^1(B)$ 
satisfying $A(1)=0$ defines a left and a right connection, every gauge 
transformation transforms connections in connections.
\\Because of this problem it seems to be necessary to introduce the
following 
definition.
\begin{de} Let $\cal P$ be locally trivial QPFB and let $\alpha_{l,r}$ be a
left(right) gauge transformation respectively. A left connection defined by
$hor_l$ is called
$\alpha_l$-covariant if $D'_l$ defined by \[ D'_l:= 
\alpha_l \circ hor_l \circ d \circ \alpha^{-1}_l \] defines the left
covariant 
derivation of a left connection.
\\A right connection defined by $hor_r$ is called $\alpha_r$-covariant if 
$D'_r$ defined by \[ D'_r := \alpha_r \circ hor_r \circ d \circ 
\alpha^{-1}_r \] defines the right covariant derivation of a right
connection.
\\A left (right) connection is called covariant, if it is 
$\alpha_{l,r}$-covariant for all gauge transformations $\alpha_{l,r}$.
\end{de}
\begin{de} Let ${\cal G}_{l,r} \subset G_{l,r}$ be a subgroup. A left
(right) connection is ${\cal G}_{l,r}$-covariant if it is
$\alpha_{l,r}$-covariant
for all $\alpha_{l,r} \in {\cal G}_{l,r}$. \end{de}
\begin{pr} Let $\cal P$ be a locally trivial QPFB and let $\Gamma({\cal
P})$
be a differential structure of $\cal P$ where the differential algebra 
$\Gamma(H)$ is the universal one. All left (right) connections are
covariant. \end{pr}
For the proof see the remarks above.
\begin{pr} Let $\cal P$ be a locally trivial QPFB and let $\Gamma({\cal
P})$ be a differential structure on $\cal P$ where $\Gamma(H)$ is bicovariant, i.e.
the corresponding right ideal $R$ is Ad-invariant. Let ${\cal Q} \subset
G_{l,r}$ be the 
subgroup of all gauge transformations which are differentiable algebra 
isomorphisms. 
All left and right connections are ${\cal Q}$-covariant. 
\end{pr}
Proof: It is sufficient to prove this assertion on a trivial bundel $B
\otimes H$.
As noted at the beginning of Section \ref{gauge}, an algebra 
automorphism $\alpha : B \otimes H \longrightarrow B \otimes H$ which is a 
gauge transformation corresponds to a homomorphism $\tau_{\alpha} : H 
\longrightarrow B$ with the property that $\tau_{\alpha}(H)$ lies in the
center of $B$. 
We have assumed that $\alpha$ is differentiable with respect to $\Gamma(B) 
\hat{\otimes} \Gamma(H)$. Let $J(B \otimes H) \subset \Omega(B \otimes H)$
be 
the differential ideal corresponding to $\Gamma(B) \hat{\otimes}
\Gamma(H)$. The 
assumption that $\alpha$ is differentiable means 
$\alpha_{\Omega}(J(B \otimes H))=J(B \otimes H)$. As shown in \cite{cama1},
Proposition 4, $J(B
\otimes H)$ is generated by the sets \begin{eqnarray*} 
&& (id \otimes 1)_{\Omega}(J(B));~ \{\sum (1 \otimes S^{-1}(r_2))d(1 
\otimes r_1)|\; r \in R \};\\
&& \{ (a \otimes 1)d(1 \otimes h)-(d(1 \otimes h))(a \otimes 1)|
\; a \in B,\,h \in H \}. \end{eqnarray*}
Applying $\alpha_{\Omega}$ to these sets and using the Ad-invariance of $R$
gives the following identities in $\Gamma(B)$.
\begin{eqnarray*}\sum \tau_{\alpha^{-1}}(r_1)d\tau_{\alpha}(r_2)&=&0, \;
\forall r \in R \\
(da) \tau_{\alpha}(h)&=& \tau_{\alpha}(h) da , \;\forall a \in B, \, \forall
h \in H. \end{eqnarray*}
Let $A_{l,r}$ be the linear map corresponding to a connection on $B \otimes
H$. Using formula (\ref{transfa}) one obtains for the transformed linear 
maps $A'_{l,r}$ 
\[ A'_{l,r}(h) = \sum A_{l,r}(h_2) \tau_{\alpha}(S(h_1)h_3) + 
\sum \tau_{\alpha^{-1}}(h_1)d\tau_{\alpha}(h_2). \] 
Inserting an element $r \in R$
in this equation and using the Ad-invariance of R, i.e. $\sum r_2 \otimes 
S(r_1)r_3 \in R \otimes H$, and the properties of $\tau_{\alpha}$ one
obtains \[A'_{l,r}(r)=0,\; \forall r \in R.\] Thus, $A'_{l,r}$ defines a left
(right) connection again. \hfill$\Boxneu$
\begin{pr}  \label{dcov} Let $\cal P$ be a locally trivial QPFB. Assume there is given a 
differential structure on  $\cal P$. Let $D_{l,r}$ be a left (right) covariant 
derivative and $\alpha_{l,r}$ a left (right) gauge transformation. The Map 
$D'_{l,r}$ defined by \[ D'_{l,r}:=\alpha_{l,r} \circ 
D_{l,r} \circ \alpha^{-1}_{l,r} \] is a left (right) covariant derivative. 
\end{pr} We leave the proof for the reader.
\begin{pr} Let $\cal P$ be a locally trivial QPFB, let $D_{l,r}$ be a 
covariant derivative and let $\alpha_{l,r}$ be a gauge transformation. \\
Then there are the following 
transformation formulas for the left (right) curvature form $\Omega_{l,r}$:
\begin{eqnarray}  {\Omega'}_l(h) &=& \sum g_{\alpha^{-1}_l}(h_1)
\Omega(h_2) 
g_{\alpha_l}(h_3) \\ {\Omega'}_r(h) &=& \sum g_{\alpha_r}(h_3) \Omega(h_2) 
g_{\alpha^{-1}_r}(h_1) \end{eqnarray} \end{pr}
Proof: Using formulas (105) and (106) of \cite{cama1} and (\ref{curv3}) and 
(\ref{curv4}) one obtains 
\[ (D'_{l,r})^2 = \alpha_{l,r} \circ (D_{l,r})^2 
\alpha^{-1}_{l,r} \] and 
\begin{eqnarray*}  \alpha_l \circ ( D_l)^2 \circ \alpha^{-1}_l
(f) &=&
\sum f_0 g_{\alpha^{-1}_l}(f_1) \Omega_l(f_2) g_{\alpha_l}(f_3) \\
\alpha_r \circ (D_r)^2 \circ \alpha^{-1}_r(f) &=& 
\sum g_{\alpha_r}(f_3) \Omega_r(f_2) g_{\alpha^{-1}_r}(f_1)
f_0.\end{eqnarray*}
for $f \in {\cal P}$. \hfill$\Boxneu$
\\At the end of this section let us remark that on locally trivial QVB
every gauge transformation $\eta$ transforms connections $\nabla$ in 
connections $\nabla'$ 
by 
\[ \nabla'= \eta \circ \nabla \circ \eta^{-1}. \]  
This is analogous to 
Proposition \ref{dcov}.
\section{Example}
This example is constructed to show that there exist nonclassical
gauge transformations. 
\\Here we glue together a noncommutative ``tube'' with the quantum group
$SU_{\nu}(2)$ along the classical subspace $S^1$ and construct a 
$SU_{\nu}(2)$ bundle over this ``base''. 
First, we define the algebra over the noncommutative ``tube'' as the
algebra $B_1$ generated by the elements
 \[x,x^*,y=y^*\] satisfying the relations
\begin{eqnarray*} xx^* &=& x^*x=1\\
 xy &=& qyx  \\   x^*y &=& q^{-1} y x^*, \end{eqnarray*} 
  where $q \in (0,1]$. It is easy to see that there exists a surjective 
  homomorphism $\pi^1_2 : B \longrightarrow P(S^1)$ defined by
  \begin{eqnarray*}
  \pi^1_2(x) &=& a\\
  \pi^1_2(x^*) &=& a^* \\
  \pi^1_2(y) &=& 0, \end{eqnarray*}
where $a$ is the generator of $P(U(1))$.
\\There also exists a surjective homomorphism
$\pi^2_1: 
  P(SU_{\nu}(2)) \longrightarrow P(S^1) $ defined by 
  \begin{eqnarray*} \pi^2_1(\alpha) &=& a \\
  \pi^2_1(\alpha^*) &=& a^* \\
  \pi^2_1(\gamma) &=& \pi^2_1(\gamma^*)=0, \end{eqnarray*}
where $\alpha,\gamma$ are the usual generators of $P(SU_{\nu}(2))$.
Our basis algebra $B$ is defined as the gluing of $B_1$ and $P(SU_{\nu}(2))$ by means of $\pi^1_2$ and $\pi^2_1$,
  \[B:=\{ (f_1,f_2) \in B_1 \bigoplus P(SU_{\nu}(2))| \pi^1_2(f_1)
= \pi^2_1(f_2) \}. \]
  Now one chooses the transition functions $\tau_{ij} : P(SU_{\nu}(2) 
  \longrightarrow P(S^1)$ as follows:
 \begin{eqnarray*}
  \tau_{12}(\alpha) &=& a \\
  \tau_{12}(\alpha^*) &=& a^* \\
  \tau_{12}(\gamma) &=& \tau_{12}(\gamma^*) = 0 \end{eqnarray*}
  and obtains a locally trivial QPFB $\cal P$ with  ``structure group''
  $SU_{\nu}(2)$.
  \\According to Proposition \ref{gauone}, to construct a left gauge
transformation 
  $\alpha_l$ we have to find a linear map $g_{\alpha_l}: P(SU_{\nu}(2)) 
  \longrightarrow {\cal P}$ fulfilling (\ref{gauone1}) -
(\ref{gauone5}.
\\First we define the linear maps $\tau^{(n)}_1 : P(SU_{\nu}(2)) 
  \longrightarrow B_1$ and $\tau^{(n)}_2: P(SU_{\nu}(2)) \longrightarrow
  P(SU_{\nu}(2)) $, $n\geq 1$. $\tau^{(n)}_1$  is assumed to be a homomorphism
defined by 
  \begin{eqnarray*}
  \tau^{(n)}_1(\alpha) &=& x^n \\
  \tau^{(n)}_1(\alpha^*) &=& {x^*}^n \\
  \tau^{(n)}_1(\gamma)&=& \tau^n_1(\gamma^*)=0. \end{eqnarray*}
  $\tau^{(n)}_2$ is defined as follows:
  \[ \tau^{(n)}_2(h)= \sum h_1h_2...h_n. \]
  The linear maps $\tau^{(n)}_i$ are convolution invertible with convolution
inverse 
  \[{\tau^{(n)}_1}^{-1}(h)=\tau^{(n)}_1(S(h)) \] and 
  \[{\tau^{(n)}_2}^{-1}(h)=\sum S(h_1)S(h_2)...S(h_n). \]
  By an easy calculation one obtains the identity 
  \[ \pi^1_2(\tau^{(n)}_1(h)) = \sum \tau_{12}(h_1) \pi^2_1(\tau^{(n)}_2(h_2)) 
  \tau_{21}(h_3)=\pi^2_1(\tau^{(n)}_2(h)),
\] 
hence (see formula (\ref{gei}) )
there exists a convolution invertible map $g^{(n)}_{\alpha_l} : P(SU_{\nu}(2))
\longrightarrow {\cal P}$ which determines a left gauge transformation.


\begin{thebibliography}{100}
\bibitem{brma} Brzezi\'nski, T., and S. Majid: Quantum group gauge theory
on quantum spaces, {\it Commun. Math. Phys.} {\bf 157} (1993), 591-638,
{\sf hep-th/9208007},
Preprint DAMTP/92-27
\bibitem{Kon} Budzy\'nski, R. J. and W. Kondracki: Quantum principal fiber
bundles: Topological aspects, {\it Rep. Math. Phys.} {\bf 37} (1996),
365-385, preprint 517 PAN Warsaw 1993, {\sf hep-th/9401019}
\bibitem{cama} Calow, D. and R. Matthes: Covering and gluing of algebras and
differential algebras, {\it J. Geom. and Phys.} {\bf 32} (2000), 364--396, {\sf
math.QA/9910031}, Preprint NTZ 25/1998
\bibitem{cama1} Calow, D. and R.
Matthes: Connections on locally trivial  quantum principal fibre bundles,
submitted to {\it J. Geom. and Phys.}, {\sf math.QA/0002228}
\bibitem{cama2} Calow, D. and R. Matthes: Locally trivial quantum vector
bundles
and associated vector bundles, {\sf math.QA/0002229}
\bibitem{dipl} Calow, D.: Differentialkalk\"ule auf Quantengruppen,
Diplomarbeit,
Leipzig 1995
\bibitem{haj} Hajac, P. M.: Strong connections on quantum principal
bundles, {\it Commun. Math. Phys.} {\bf 182} (1996), 579--617
\bibitem{korn} Koornwinder, T. H.: General Compact Quantum Groups, a
Tutorial,
University of Amsterdam, Faculty of Mathematics and Computer Science
\bibitem{wor3} Woronowicz, S. L.: Compact matrix pseudogroups,
{\it Commun.
Math. Phys.} {\bf 111} (1987), 613--665
\bibitem{wor2} Woronowicz, S. L.: Differential calculus on compact
matrix pseudogroups (quantum groups), {\it Commun. Math. Phys.} {\bf 122}
(1989), 125--170
\end{thebibliography}
\end{document}